\documentclass[12 pt, draftcls, onecolumn]{IEEEtran}
\usepackage{pdfsync}
\usepackage{amssymb,epsfig,color,cite,amsmath,amsfonts,mathrsfs,algorithm,algorithmic}
\usepackage{epstopdf}

\newtheorem{lemma}{Lemma}[section]
\newtheorem{theorem}{Theorem}[section]

\newtheorem{remark}{Remark}[section]

\newtheorem{assumption}{Assumption}[section]
\allowdisplaybreaks

\begin{document}
\title{Dynamical Primal-Dual Nesterov Accelerated Method and Its Application to Network Optimization}
\title{%\huge
Dynamical Primal-Dual Accelerated  Method  with Applications to Network Optimization
}

\author{Xianlin~Zeng, Jinlong Lei,  and Jie~Chen %<-this % stops a space
%\thanks{This work was supported by National Key Research and Development Program of China under Grant 2018AAA0101601,  Shanghai Municipal Commission of Science and Technology Project No. 19511132101, NSFC  (Grant no.  62073035, 61873033, 61991414, 62003240),  Projects of Major International (Regional) Joint Research Program NSFC (Grant no. 61720106011),   Research Fund of Science and Technology on Space Intelligent Control Laboratory (Grant no. 6142208200312), and Shanghai Sailing Program.}
 \thanks{X. Zeng (xianlin.zeng@bit.edu.cn) is with the Key Laboratory of Intelligent Control and Decision of Complex Systems, School of Automation, Beijing Institute of Technology, 100081, Beijing, China.}
 \thanks{J. Lei  (leijinlong@tongji.edu.cn) and J. Chen (chenjie@bit.edu.cn) are with the Department of Control Science and Engineering, Tongji University, Shanghai, China;  the Shanghai Research Institute for Intelligent Autonomous Systems, Shanghai, China. J. Chen is also with the Key Laboratory of Intelligent Control and Decision of Complex Systems, School of Automation, Beijing Institute of Technology, 100081, Beijing,  China}
}

%\markboth{Submitted to IEEE Transactions on Automatic Control}{}

\maketitle

\begin{abstract}
This paper develops a  continuous-time  primal-dual accelerated method with  an increasing damping coefficient for a class of convex optimization problems with affine equality constraints. % and applies the method to two classes of network optimization problems.
% \blue{To be specific},   a \blue{modified version of the continuous-time accelerated  primal-dual method is designed} by using an increasing damping coefficient. %As far as we know,  this is the first continuous-time Nesterov  accelerated primal-dual method.
  This paper analyzes critical values for parameters in the proposed  method and  prove  that the  rate of  convergence  in terms of the duality gap function is $O(\tfrac{1}{t^2})$ by choosing suitable parameters.  As far as we know,  this is the first continuous-time primal-dual   accelerated  method that can obtain the optimal   rate. Then  this work applies  the proposed   method  to  two network optimization  problems, a  distributed optimization problem with consensus constraints and   a distributed extended monotropic optimization problem, and  obtains two variant distributed algorithms.  Finally, numerical simulations are given to   demonstrate  the efficacy of the proposed method.
\end{abstract}

\begin{IEEEkeywords}
Nesterov's accelerated method, primal-dual method,  network optimization,  continuous-time algorithm.
\end{IEEEkeywords}

% \fi

\section{Introduction}

Many large-scale problems in control and optimization (including optimal consensus of multi-agents \cite{SJ:2012}, network flow \cite{ZYHX:SIAM:2018}, energy dispatch of power grids \cite{YHL:Automatica:2016}, etc.) can be formulated as distributed optimization problems  that studied in  \cite{NOP:2010,WE:2011}.
In distributed optimization, information is allocated to multiple agents with no central agent, and certain equality constraints such as consensus constraints \cite{GC:2014}, resource allocation constraints \cite{YHL:Automatica:2016}, and monotropic constraints \cite{ZYHX:SIAM:2018} should be satisfied.  Distributed optimization  problems are often solved by first-order methods, whose implementation is simpler than that of higher-order algorithms.
 The discrete-time  distributed algorithms such as distributed subgradient methods \cite{NOP:2010}, distributed Nesterov's gradient methods \cite{JXM:tac:2014},  distributed gradient tracking methods \cite{pu2018push}, and distributed  primal-dual   methods  \cite{YHL:SCL:2015,lei2016primal}  have been extensively studied, of which  the convergence analysis is specific to that  designed algorithm while not within a unified framework. Recently, the study of the  continuous-time  distributed optimization algorithms  have also drawn much  attention from researchers, % by constructing Lyapunov functions  for the convergence analysis,
 see e.g.,  \cite{ZYH:2017,GC:2014,YHX:SCL:2014}.

It is well-known  that  the fastest rate of convergence of first-order methods for convex optimization in the worst case  is  $O(\tfrac{1}{t^2})$, see \cite{Nesterov:intro:2004}. Most existing approaches with  rate $O(\tfrac{1}{t^2})$ mainly focus on primal algorithms, which, however,  can not be directly  applied to the primal-dual framework.  Because the primal-dual framework is a  powerful technique for constrained optimization with equality/inequality constraints, which primal-based methods can not deal with, accelerating  the rate of convergence of primal-dual methods is of great importance.
The primary objective of this paper is to propose  a dynamical primal-dual   accelerated method for convex network optimization with convergence rate  $O(\tfrac{1}{t^2})$. The previous works  demonstrated that primal-dual methods for distributed optimization   can guarantee  asymptotic convergence (see \cite{ZYH:2017}) or a rate of convergence $O(\tfrac{1}{t})$  (see \cite{ZYHX:SIAM:2018}). %On the other hand,  Nesterov accelerated method and its variants, which have a rate of convergence $O(\frac{1}{t^2})$, mainly focus on optimization problems without equality constraints.
Hence, for modern large-scale network optimization problems with   equality constraints, the   convergence rate  of existing distributed primal-dual methods is slower than that of algorithms for centralized  convex optimization  without equality constraints. Thus, it is important to design a primal-dual accelerated method for   convex network  optimization problems.

\subsection{Related Work}

\subsubsection{Distributed optimization algorithms}
Due to the  rapid growth in the scale and complexity of  the network  optimization problems, distributed first-order  gradient-based primal-dual methods have been gaining more attention because they are easy for distributed implementation.
For distributed strongly convex optimization problems, distributed solvers with linear or exponential rates   have been proposed in \cite{KCM:A:2015,YHL:Automatica:2016}. Furthermore, by using the notion of metric subregularity,   some existing distributed algorithms (see \cite{LWY:tac:2019}) have been proved  to have linear convergence rates if non-strongly convex cost functions satisfy some properties.  If cost functions are arbitrary convex functions, distributed primal-dual algorithms are only proved to have asymptotic convergence  (see \cite{YHL:SCL:2015,ZYH:2017,YWL:tac:2019,LYH:TAC:2017,Kia:SCL:2017,MK:arxiv:2020}) or convergence   rate  $O(\frac{1}{t})$ (see \cite{ZYHX:SIAM:2018}). Recently, a distributed algorithm that achieves convergence rate $O(\frac{1}{t^{1.4-\epsilon}})$ has been proposed in \cite{JXM:tac:2014}, and  an accelerated distributed  algorithm using a gradient estimation scheme with rate $O(\frac{1}{t^{2-\epsilon}})$  has been proposed in \cite{QL:Arxiv:2019}. However, both \cite{JXM:tac:2014} and \cite{QL:Arxiv:2019}  consider consensus constraints, cannot achieve  the optimal rate  $O(\tfrac{1}{t^2 })$,  and are not primal-dual  algorithms.

% P1 Nesterov Accelerated method
\subsubsection{Centralized Nesterov's accelerated gradient methods} The Nesterov accelerated method,   using a vanishing damping coefficient,  was  developed  in  \cite{Nesterov} and proved to have a rate of convergence $O(\frac{1}{t^2 })$.
As a first-order algorithm, the  Nesterov's accelerated method and its different variants (see \cite{Nesterov:intro:2004}) are proved to be optimal in some sense (see \cite{NY:1983}) and have been widely studied in different settings \cite{BX:SIAM:2009,Lan:MP:2012}. Recently, growing attention has been dedicated to the design and analysis of continuous-time Nesterov's accelerated methods (see \cite{Su2015A,WWJ:2016,ACR:ESAIM:2017,SIEGEL:Arxiv:2019}). On one hand, ordinary differential equations (ODEs) often exhibit similar convergence properties to their discrete-time counterparts  and thus can serve as a tool for algorithm design and analysis. On the other hand, continuous-time algorithms may allow for a better understanding of  intuitive and ideas in the design. \cite{Su2015A} showed that the continuous-time counterpart  of the Nesterov's accelerated method is a second-order ODE, and  proved that the suboptimality gap is $O(\frac{1}{t^2})$ for parameter $\alpha\geq 3$ in the ODE;  \cite{ACPR:MP:2018} further proved that  the generated trajectory   converges to a minimizer of $\phi$ as $t\rightarrow\infty$ for $\alpha>3$; while  \cite{ACR:ESAIM:2017} proved that the rate of convergence of the ODE  proposed in \cite{Su2015A} is $O({t^{-\frac{2\alpha}{3}}})$ for $0<\alpha\leq 3$.

To develop new insights into accelerated algorithms and obtain new algorithms, there are ongoing recent research on discretization of accelerated algorithms.
\cite{KBB:NIPS:2015} studied a mixed forward/backward Euler scheme  for continuous-time algorithms for constrained optimization and proved an analogous $O(\frac{1}{k^2})$ rate when the step size is small enough.
%We show if a Hessian damping term is added to the ODE from Su et al. (2014), then Nesterov¡¯s method arises as a straightforward discretization of the modifified ODE.
\cite{SDJS:arxiv:2018} proposed  alternative ODEs of Nesterov's accelerated methods and their discretizations, which are called high-resolution ODEs because they use Hessians of cost functions to distinguish between Nesterov's accelerated gradient method for strongly convex functions and Polyak's heavy-ball method.
\cite{SDSJ:NIPS:2019} considered three discretization schemes: an explicit Euler scheme, an implicit Euler scheme,
and a symplectic scheme, and applied the symplectic scheme to a high-resolution ODE proposed by \cite{SDJS:arxiv:2018} for minimizing smooth strongly convex functions.
\cite{LO:arxiv:2020} showed that the Nesterov's method arises as a straightforward discretization of a modified ODE. \cite{SGK:Arxiv:2020} further derived an ODE model of the accelerated triple momentum algorithm for   strongly convex optimization  and investigated the convergence behavior of the ODE model.
However, these works only focus on unconstrained optimization problems, while can not be applied to the widely used primal-dual framework.

In recent years, accelerated primal-dual methods have received lots of attention due to its application in constrained optimization and minimax optimization (or saddle point problems). Focusing on    duality gaps of the Lagrangian function or minimax optimization, which often converge to 0 with rate $O(\frac{1}{k})$,  various accelerated primal-dual algorithms   have been investigated.
In \cite{Chambolle:2011:FPA:1968993.1969036}, a primal-dual algorithm with a rate of convergence $O(\frac{1}{k^2})$ has been proposed for saddle point problems, whose primal or dual cost function is uniformly convex. In \cite{Chambolle:2011:FPA:1968993.1969036}, since each iteration involves solving an optimization problem,   cost functions need to be ``simple". %\cite{CLO:SIAM:2014} proposed an accelerated primal-dual method for solving a class of deterministic and stochastic saddle point problems by incorporating a multistep acceleration scheme  without smoothing the objective function.
\cite{XZ:COA:2018} proposed an accelerated design for primal-dual block coordinate method that has a rate of convergence  $O(\frac{1}{k^2})$ for strongly convex functions.
\cite{YXu:SIAM:2018} proposed a linearized augmented Lagrangian method with full acceleration $O(\frac{1}{k^2})$ for strongly convex functions with adaptive parameters.

Different from previous literature, this work is the {\em first attempt} to propose a dynamical primal-dual  method for constrained convex optimization problems  with  a rate of convergence $O(\frac{1}{t^2})$.

\subsection{Contributions}
%This paper aims to propose a primal-dual Nesterov accelerated method for convex optimization problems with equality constraints and apply the method to  distributed optimization problems.
This paper has  made the following  contributions.
\begin{enumerate}
  \item Considering convex optimization problems with equality constraints, this paper proposes a primal-dual Nesterov's accelerated method, which has  a convergence rate $O(\frac{1}{t^2})$.
      The novel part in the design is the use of the derivative information in the Lagrangian saddle point dynamics to obtain the accelerated convergence of the proposed first-order primal-dual method. To our best knowledge, this is the first continuous-time primal-dual accelerated method owning  a convergence rate $O(\frac{1}{t^2})$.
     % Different from  primal-based accelerated distributed algorithms  \cite{JXM:tac:2014,QL:Arxiv:2019}, this method can be applied to problems with more general linear equality constraints.

  \item This paper further analyzes convergence properties of   the proposed primal-dual   accelerated method for difference choices of parameters, and figures out  the best choice of parameters  for the optimal convergence rate $O(\frac{1}{t^2})$.
       To be specific, $\alpha>3$ and $\beta = \frac{1}{2}$ are shown in Section \ref{method}.%Both the algorithm design and convergence properties are consistent with existing accelerated algorithms.

  \item This paper applies the proposed primal-dual accelerated method to two classes of widely studied network optimization problems: distributed optimization with consensus constraints and distributed extended monotropic optimization. This
       leads to two distributed primal-dual Nesterov's accelerated algorithms with convergence guarantees. The numerical experiments  show faster convergence performances  than that of existing results on distributed optimization problems \cite{WE:2011,ZYHX:SIAM:2018}.
\end{enumerate}

\subsection{Organization}

This paper is organized as follows. Section \ref{sec:background} provides necessary mathematical preliminaries and formulates  a convex optimization problem with linear  equality constraints. Section \ref{method} proposes a continuous-time primal-dual Nesterov's accelerated method
  and gives  its convergence properties. Based on the proposed method, Section \ref{DO}
further designs two primal-dual Nesterov's accelerated algorithms for   two widely studied  network optimization problems  and presents the
simulation results. Section \ref{conclusion}  concludes this paper.

\section{Mathematical Preliminary and Problem Setup}\label{sec:background}
In this section, we introduce some mathematical notations and give the problem statement.
\subsection{Notation}
\label{sec:def}
The symbol  $\mathbb R$ denotes the set of real numbers; %$\overline{\mathbb R}_+$ denotes the set of nonnegative real numbers,
$\mathbb R^n$ denotes the set of $n$-dimensional real column vectors; $\mathbb R^{n\times m}$ denotes the set of $n$-by-$m$ real matrices; $I_n$ denotes the $n\times n$ identity matrix;  $(\cdot)^\mathrm{T}$ denotes transpose.
 We write $\mathrm {rank}\,A$ for the rank of the matrix $A$,  $\mathrm{range} (A)$ for the range of the matrix $A$,
 $\ker (A)$ for the kernel of the matrix $A$, %$\lambda _{\mathrm {max}} (A)$ for the largest eigenvalue of the matrix $A$,
 $\mathbf 1_n$ for the $n\times 1$ ones vector, $\mathbf 0_n$ for the $n\times 1$ zeros vector, and $A \otimes B$ for the Kronecker product of matrices $A$ and $B$.
Furthermore, $\|\cdot\|$ denotes the Euclidean norm; % $\|\cdot\|_p$ denotes the $p$-norm where $p\geq 1$;
 $A> 0$ ($A\geq 0$) denotes that matrix $A\in\mathbb R^{n\times n}$ is positive definite (positive semi-definite);
%$\|x(t)\|_{\infty} =\max_{t\geq 0}\|x(t)\|^2$, $\overline{\mathcal{S}}$ denotes the closure of the subset $\mathcal{S}\subset\mathbb{R}^{n}$; $\mathrm{int}(\mathcal{S})$ denotes the interior of the subset $\mathcal{S}$; %$\mathrm{dim}(\mathcal S)$ denotes the dimension of the vector space $\mathcal S$;
%$\mathcal B_{\epsilon}(\alpha),\alpha\in\mathbb R^n,\epsilon>0,$ denotes the open ball {\em centered} at $\alpha$ with {\em radius} $\epsilon$; $\mathrm{dist}(p,\mathcal M)$ denotes the distance from the point $p$ to the set $\mathcal M$, that is, $\mathrm{dist}(p,\mathcal M) \triangleq \inf_{x\in\mathcal M}\|p-x\|$; $x(t)\rightarrow \mathcal M$ as $t\rightarrow \infty$ denotes that $x(t)$ approaches the set $\mathcal M$, that is, for each $\epsilon>0$ there exists $T>0$ such that $\mathrm{dist}(x(t),\mathcal M)<\epsilon$ for all $t>T$.
 Let $f:\overline{\mathbb R}_+\rightarrow\overline{\mathbb R}_+$ be a continuous  function. $f(t) = O(\frac{1}{t^n})$ indicates that there exist constants $C>0$ and $t_0\geq 0$ such that $f(t)\leq \frac{C}{t^n}$ for all $t\geq t_0$.

A set $\Omega $ is {\em convex} if $\lambda z_1
+(1-\lambda)z_2\in \Omega$ for all $z_1, z_2 \in \Omega$ and $\lambda\in [0,\,1]$.  A function $f: \Omega \to \mathbb{R}$ is  {\em convex (strictly convex)}  if $f(\lambda z_1
+(1-\lambda)z_2) \leq (<) \lambda f(z_1) + (1-\lambda)f(z_2)$ for any $z_1,z_2\in \Omega$, $z_1\not= z_2$,  and $\lambda\in (0,\,1)$.
If function $f: \Omega \to \mathbb{R}$ is (strictly) convex, it is well-known that $(z_1-z_2)^{\top}(\nabla f(z_1)-\nabla f(z_2))\geq (>) 0$ for any $z_1,z_2\in \Omega$ and $z_1\not= z_2$. Given a differentiable function $f(x,y)$,  $\nabla_{x}f(x,y)$ denotes the partial gradient of function $f(x,y)$ with respect to $x$.

An undirected graph $\mathcal G$ is denoted by $\mathcal G(\mathcal V, \mathcal E, A)$, where $\mathcal V=\{1,\ldots, n\}$ is a set of nodes, $\mathcal E\subset\mathcal V \times \mathcal V$ is a set of edges, and $ A=[a_{i,j}]\in\mathbb R^{n\times n}$ is an {\em adjacency matrix} such that $a_{i,j}=a_{j,i}>0$ if $(j,i)\in\mathcal E$ and $a_{i,j}=0$ otherwise.
The {\em Laplacian matrix} is $L_n=D- A$, where $D\in\mathbb R^{n\times n}$ is diagonal with $D_{i,i}=\sum_{j=1}^n a_{i,j}$, $i\in\{1,\ldots,n\}$.
If the graph $\mathcal G$ is undirected and connected,
then $L_n=L_n^{\top}\geq 0$, $\mathrm{rank} \,L_n=n-1$ and $\ker (L_n)=\{k\mathbf{1}_n:k\in\mathbb R\}$.

\subsection{Problem Formulation}

Consider a convex optimization problem with an affine equality constraint given by
\begin{eqnarray}\label{CP}
\min_{x\in\mathbb R^q}\ \phi(x),\quad \text{ s.t. }Ax-b=0_m,
\end{eqnarray}
where $A\in\mathbb R^{m\times q}$, $b\in\mathbb R^m$, and $\phi:\mathbb R^q\rightarrow\mathbb R$ is a convex and twice differentiable cost function.

It follows from the Karush-Kuhn-Tucker (KKT) optimality condition (see \cite [Theorem 3.34]{Ruszczynski:2006}) that $x^*\in\mathbb R^q$ is a solution to the problem \eqref{CP} if and only if there exists $\lambda^*\in\mathbb R^m$ such that
\begin{subequations}\label{KKT}
	\begin{align}
		\mathbf{0}_q  = & \nabla \phi(x^*)+A^{\top}\lambda^*,\label{KKT1}\\
		\mathbf{0}_m =& Ax^*-b.
	\end{align}
\end{subequations}
Define the augmented Lagrangian function $L:\mathbb R^q\times \mathbb R^m\rightarrow\mathbb R$ as
\begin{eqnarray}\label{Lag1}
	L(x,\lambda) = \phi(x)+\lambda^{\top}(Ax-b)+\frac{1}{2}\|Ax-b\|^2.
\end{eqnarray}
It is well-known that $(x^*,\lambda^*)\in\mathbb R^q\times \mathbb R^m$ satisfies \eqref{KKT} if and only if $(x^*,\lambda^*)\in\mathbb R^q\times \mathbb R^m$ is a saddle point of \eqref{Lag1}, that is,
$$L(x,\lambda^*)\geq L(x^*,\lambda^*)\geq L(x^*,\lambda),\quad \forall (x,\lambda)\in\mathbb R^q\times \mathbb R^m.$$

Given a continuous-time algorithm,  the rate of convergence of the algorithm is said to be $O(\frac{1}{t^p})$ if the duality gap satisfies
$L(x(t),\lambda^*)-L(x^*,\lambda^*)=O(\frac{1}{t^p})$, where $p>0$.
This paper aims to design a primal-dual accelerated method that has a convergence rate faster than  $O(\frac{1}{t})$, the rate of first-order primal-dual methods for convex optimization.

%\begin{remark}
%Problem \eqref{CP}   has seen wide applications in  machine learning,
%data mining, and image processing. The primal-dual framework is known to be efficient for constrained optimization problems like \eqref{CP}.  In recent years, with the rise of big data problems, first-order methods have received tremendous attention because they are computationally cheap and easily to be implemented in a distributed or parallel way. However, the first-order primal-dual methods  converge at a slow rate  $O(1/t)$. Hence, the development of primal-dual accelerated algorithms is of great importance.\hfill $\Diamond$
%\end{remark}

Furthermore, the formulation  \eqref{CP} captures   two important scenarios of network optimization problems.

{\em Scenario 1: Distributed Optimization with Consensus Constraints. }Consider a network of $n$ agents interacting over a graph $\mathcal G$.
The distributed agents cooperate  to solve the following problem
\begin{subequations}\label{S1}
\begin{align}
\min_{x\in\mathbb R^{nq}}&f(x), \quad f(x)= \sum_{i=1}^n f_i(x_i),\\
\mathrm{s.t.}~&x_i=x_j,\quad i,j\in\{1,\ldots,n\},
\end{align}
\end{subequations}
where agent $i$ only knows  its   local cost function  $f_i:\mathbb{R}^q\to \mathbb{R}$  and the shared information of its neighbors through  local communications.

Problem \eqref{S1} is a widely investigated model that has  many applications such as optimal consensus of agents \cite{SJ:2012} %, routing  of wireless sensor networks \cite{ML:2006},
  and distributed machine learning  \cite{BPCPE:2011}.

{\em Scenario 2: Distributed Extended Monotropic Optimization. } Given a network $\mathcal G$ composed of $n$ agents, the distributed extended monotropic optimization problem is
\begin{subequations}\label{S2}
\begin{align}
\min_{y\in\mathbb R^q}&\,{h}(y), \quad h(y)= \sum_{i=1}^n {h}_i(y_i),\label{S2-1}\\
\mathrm{s.t.}~ & Wy=\sum_{i=1}^n W_i y_i =\sum_{i=1}^n d_i= d_0,\quad i\in\{1,\ldots,n\},\label{optimization_problem2}
\end{align}
\end{subequations}
where  $y_i\in\mathbb R^{q_i}$,  $y \triangleq [y_1^{\top},\ldots,y_n^{\top}]^{\top}\in\mathbb R^{q}$ with $q\triangleq \sum_{i=1}^{n}q_i $, $W_i\in\mathbb{R}^{m\times q_i}, $ $d_0,d_i\in\mathbb R^{m}$, and $W = [W_1,\ldots,W_n]\in\mathbb{R}^{m\times q}$. In this problem, agent $i$ has its state $y_i\in\mathbb R^{q_i}$, objective function $h_i(y_i)$,  constraint matrix $W_i\in\mathbb{R}^{m\times q_i}$,  information from neighboring agents, and  a vector $d_i\in\mathbb R^m$ such that $\sum_{i=1}^n d_i = d_0$.

Problem (\ref{S2}) covers many network optimization problems such as resource allocation problems and network flow problems \cite{YHL:SCL:2015,YHL:Automatica:2016,ZYHX:SIAM:2018}. Different from problem (\ref{S1}),  $y_i$'s appear in the same constraint \eqref{optimization_problem2}.

To ensure the wellposedness of  problems \eqref{S1} and \eqref{S2}, the following assumption is needed.
\begin{assumption}\label{A1}
\begin{enumerate}
  \item Graph $\mathcal G$ is connected and undirected.
  \item There exists at least one finite solution to problems (\ref{S1}) and \eqref{S2}.
  \item  Function $f_i(\cdot)$ in problem (\ref{S1}) ($h_i(\cdot)$ in problem \eqref{S2}) is convex  and twice continuously differentiable.
\end{enumerate}
\end{assumption}

\section{Primal-Dual Nesterov's Accelerated Method}\label{method}
In this section, we propose a dynamical primal-dual Nesterov's accelerated  method for problem \eqref{CP}.
We then investigate convergence properties of the proposed method and analyze  critical values of algorithm  parameters.
\subsection{Algorithm Design}

%Let $(x^*,\lambda^*)\in\mathbb R^q\times \mathbb R^m$.

%
%An intuitive  extension of the primal-based Nesterov accelerated  algorithm \eqref{AGV} to a primal-dual case for problem \eqref{CP} may be
%\begin{subequations}\label{natural1}
%\begin{align}
%\ddot x(t)= & -\frac{\alpha}{t}\dot x(t)-\nabla_{x(t)} L(x(t),\lambda(t)) , \, x(t_0)=x_0,\quad \dot x(t_0)=\dot x_0,\\
%\dot \lambda(t)=&  \nabla_{\lambda(t)} L(x(t),\lambda(t)) ,\quad \lambda(t_0)=\lambda_0,
%\end{align}
%\end{subequations}
%or
%\begin{subequations}\label{natural2}
%\begin{align}
%\ddot x(t)= & -\frac{\alpha}{t}\dot x(t)-\nabla_{x(t)} L(x(t),\lambda(t)), \quad x(t_0)=x_0,\, \dot x(t_0)=\dot x_0,\\
%\ddot \lambda(t)=& -\frac{\alpha}{t}\dot \lambda (t)+ \nabla_{\lambda(t)} L(x(t),\lambda(t)),\quad \lambda(t_0)=\lambda_0,\, \dot \lambda(t_0)=\dot \lambda_0.
%\end{align}
%\end{subequations}
%However, the involvement of the dual variable $\lambda(t)$ makes the analysis of completely different from that of primal-based Nesterov accelerated  algorithm \eqref{AGV}.  Those intuitive extensions may not be right primal-dual designs for the Nesterov  accelerated algorithm. To our best knowledge, there are no theoretical results for algorithms \eqref{natural1} or \eqref{natural2}.
We propose a primal-dual accelerated method as follows:
\begin{subequations}\label{PD-Cent}
\begin{align}
\ddot x(t)= & -\frac{\alpha}{t}\dot x(t)-\nabla \phi(x(t))-A^{\top}(\lambda(t)+\beta t\dot \lambda(t))-A^{\top}(Ax(t)-b), \\
\ddot \lambda(t)=& -\frac{\alpha}{t}\dot \lambda(t) + A(x(t)+\beta t\dot x(t))-b,
\end{align}
\end{subequations}
where $t_0>0$, $\alpha>3$,  $\beta=\frac{1}{2}$, $x(t_0)=x_0,$ $\dot x(t_0)=\dot x_0$, $\lambda(t_0)=\lambda_0,$ and $\dot \lambda(t_0)=\dot \lambda_0$.
In the remaining of this paper, we omit $(t)$ in the algorithm and analysis without causing confusions. For example, we use $x$ and $\dot x$ to denote $x(t)$ and $\dot x(t)$.
Then with  the definition of  $L(x,\lambda) $   in  \eqref{Lag1}, the algorithm \eqref{PD-Cent}  can be written as
\begin{eqnarray*}
\ddot x&= & -\frac{\alpha}{t}\dot x-\nabla_{x}L(x,\lambda+\beta t\dot \lambda), \quad x(t_0)=x_0,\quad \dot x(t_0)=\dot x_0,\\
\ddot \lambda&=& -\frac{\alpha}{t}\dot \lambda + \nabla_{\lambda}L(x+\beta t\dot x,\lambda),\quad \lambda(t_0)=\lambda_0,\quad \dot \lambda(t_0)=\dot \lambda_0.
\end{eqnarray*}

\begin{remark}
Since $\phi(\cdot)$ is twice  differentiable, $\nabla\phi(\cdot)$ is locally Lipschitz continuous. It follows from \cite[Theorem 2.38, pp. 96]{HC:2008} that algorithm \eqref{PD-Cent} has  a unique trajectory. The initial time $t_0>0$ avoids the singularity of the damping coefficient $\frac{\alpha}{t}$ at zero. Although algorithm \eqref{PD-Cent} uses $t\dot x(t)$ and $t\dot \lambda(t)$ in the righthand side, Section \ref{alpha>3} will show that $t\dot x(t)$ and $t\dot \lambda(t)$ are bounded for all  $t\geq t_0$. Hence,   algorithm \eqref{PD-Cent} is well defined with a bounded right-hand side.
\hfill $\Diamond$
\end{remark}

\begin{remark}
The use of derivative information $\beta t\dot x(t)$ and $\beta t\dot \lambda(t)$ in the update \eqref{PD-Cent} can be interpreted from the following perspectives. From the control perspective, it  may be viewed as a ``derivative feedback" design and plays a role as damping terms.
From the optimization perspective, $\beta t\dot x(t)$ and $\beta t\dot \lambda(t)$ point to the future moving direction of $x(t)$ and $\lambda(t)$. Thus, algorithm \eqref{PD-Cent} uses the estimated ``future" position $x(t)+\beta  t\dot x(t)$ and $\lambda(t)+\beta  t\dot \lambda(t)$.\hfill $\Diamond$
\end{remark}

\begin{remark}
Without loss of generality, the order for the rate of convergence of algorithm \eqref{PD-Cent} remains  unchanged if the right-hand side of \eqref{PD-Cent} is multiplied by any positive constant gain. However, one should be very careful choosing time-varying gains because an infinitely large increasing gain may produce unbounded variable derivatives and make  algorithms  impractical.
\hfill $\Diamond$
\end{remark}

\subsection{Convergence Analysis for $\alpha>3$ and $\beta=\frac{1}{2}$}\label{alpha>3}\label{Ageq3}
In this subsection, we investigate the  convergence results of  the algorithm \eqref{PD-Cent} when $\alpha>3$ and $\beta=\frac{1}{2}$. To be specific, the following theorem shows that algorithm  \eqref{PD-Cent} has convergence rate $O(\frac{1}{t^2})$. % when $\alpha>3$ and $\beta=\frac{1}{2}$.

\begin{theorem}\label{thm-cent}
%Denote $S$ as the set of solutions to  problem \eqref{CP} and assume $S\not=\emptyset$.
Suppose the problem \eqref{CP}  has a nonempty solution set $S$. Let $(x(t),\lambda(t))$ be a trajectory generated by the algorithm \eqref{PD-Cent} with $\alpha>3$ and $\beta = \frac{1}{2}$. Then we have the following.
\begin{itemize}
  \item [(i)] Trajectories of $(x(t),\lambda(t),t\dot x(t),t\dot \lambda(t))$ and $(\ddot x(t),\ddot \lambda(t))$ are bounded for $t\geq t_0$.
  \item [(ii)]The trajectory  $(x(t),\lambda(t),\dot x(t),\dot \lambda(t))$ satisfies the convergence properties
$L(x(t),\lambda^*)-L(x^*,\lambda^*)=O(\frac{1}{t^2})$, $\|Ax(t)-b\|^2=O(\frac{1}{t^2})$, $\|\dot x(t)\|=O(\frac{1}{t})$, and $\|\dot \lambda(t)\|=O(\frac{1}{t})$.
 \item  [(iii)] The trajectory of  $x(t)$ converges to the solution set $S$ as $t\rightarrow\infty$.
\end{itemize}
\end{theorem}
\begin{IEEEproof}  Let $(x^*,\lambda^*) \in \mathbb{R}^q \times \mathbb{R}^m$ satisfy \eqref{KKT}.
Define function
\begin{eqnarray}\label{Vc}
V(t,x,\lambda,\dot x,\dot \lambda)=V_1(t,x)+V_2(t,x,\dot x)+V_3(t,\lambda,\dot \lambda)
\end{eqnarray}
such that
\begin{eqnarray}
%L(x,\lambda) &=& \phi(x)+\lambda^{\top}(Ax-b)+\frac{1}{2}\|Ax-b\|^2,\\
V_1&=& t^2[L(x,\lambda^*)-L(x^*,\lambda^*)],\label{eq9}\\
V_2 &=& 2\|x+\beta t\dot x-x^*\|^2+2(\alpha\beta-\beta-1)\|x-x^*\|^2,\\
V_3 &=& 2\|\lambda+\beta t\dot \lambda-\lambda^*\|^2+2(\alpha\beta-\beta-1)\|\lambda-\lambda^*\|^2,
\end{eqnarray}
where $L(\cdot,\cdot)$ is defined in \eqref{Lag1}. By the property of saddle points of $L(\cdot,\cdot)$, $L(x,\lambda^*)\geq L(x^*,\lambda^*)$ for all $x\in\mathbb R^q$. Hence, function $V$ is positive definite with respect to  $(x,\lambda,t\dot x,t\dot \lambda)$ for all $t\geq t_0$.

(i) The derivatives of $V_i$'s $i=1,2,3,$ along the trajectory of  algorithm \eqref{PD-Cent} satisfy that
\begin{eqnarray}\label{cp1-DVc1}
\dot V_1
&=& 2t[\phi(x)-\phi(x^*)+\lambda^{*\top}A(x-x^*)+\frac{1}{2}\|Ax-b\|^2]\notag\\
&&+t^2[\nabla \phi(x)+A^{\top}\lambda^*+A^{\top}(Ax-b)]^{\top}\dot x,
\end{eqnarray}
\begin{eqnarray}\label{cp1-DVc2}
\dot V_2  &=& 4(x+\beta t\dot x-x^*)^{\top}((1+\beta)\dot x+\beta t\ddot x)+4(\alpha\beta-\beta-1)(x-x^*)^{\top}\dot x\nonumber\\
&=&-4\beta t(x-x^*)^{\top}(\nabla \phi(x)+A^{\top}\lambda)-4\beta t\|Ax-b\|^2-4\beta^2 t^2(x-x^*)^{\top}A^{\top}\dot \lambda\notag\\
&&+4\beta (1+\beta-\alpha\beta)t\|\dot x\|^2-4\beta^2t^2\dot x^{\top}(\nabla \phi(x)+A^{\top}\lambda)-4\beta^3t^3\dot x^{\top}A^{\top}\dot \lambda\notag\\
&&-4\beta^2t^2\dot x^{\top}A^{\top}(Ax-b),
\end{eqnarray}
\begin{eqnarray}\label{cp1-DVc3}
\dot V_3  &=& 4 (\lambda+\beta t\dot \lambda-\lambda^*)^{\top}((1+\beta)\dot \lambda+\beta t\ddot \lambda)+4(\alpha\beta-\beta-1)(\lambda-\lambda^*)^{\top}\dot \lambda\nonumber\\
&=&4\beta t(\lambda-\lambda^*)^{\top}(Ax-b)+4\beta^2t^2 (\lambda-\lambda^*)^{\top}A\dot x+4\beta(1+\beta-\alpha\beta)t\|\dot \lambda\|^2\notag\\
&&+4\beta^2t^2\dot\lambda ^{\top}(Ax-b) +4\beta^3t^3\dot x^{\top}A^{\top}\dot \lambda.
\end{eqnarray}

Plug \eqref{KKT} in \eqref{cp1-DVc2} and \eqref{cp1-DVc3}, and rearrange the terms.  We have
\begin{align}\label{p1-temp}
\dot V_2 +\dot V_3 =&-4\beta t\big((x-x^*)^{\top}(\nabla \phi(x)+A^{\top}\lambda^*)+\frac{1}{2}\|Ax-b\|^2\big)\notag\\
&-4\beta^2t^2\dot x^{\top}\big(\nabla \phi(x)+A^{\top}\lambda^*+A^{\top}(Ax-b)\big)+N
\end{align}
with $
N \triangleq 4\beta (1+\beta-\alpha\beta)t\|\dot x\|^2+4\beta(1+\beta-\alpha\beta)t\|\dot \lambda\|^2-2\beta t\|Ax-b\|^2.$

%\begin{eqnarray*}
%N &=& 4\beta (1+\beta-\alpha\beta)t\|\dot x\|^2+4\beta(1+\beta-\alpha\beta)t\|\dot \lambda\|^2\\
%&&-2\beta t\|Ax-b\|^2.
%\end{eqnarray*}

Plugging $\beta = \frac{1}{2}$ in \eqref{p1-temp}, it follows from  \eqref{cp1-DVc1} and  \eqref{p1-temp}  that
\begin{eqnarray}
\dot V& = 2t[\phi(x)-\phi(x^*)-(x-x^*)^{\top}\nabla \phi(x)]+N,\label{p1-DV-V}
\end{eqnarray}
where $N= (3-\alpha)t\|\dot x\|^2+(3-\alpha)t\|\dot \lambda\|^2- t\|Ax-b\|^2.$
%\begin{eqnarray}
%\dot V& = 2t[\phi(x)-\phi(x^*)-(x-x^*)^{\top}\nabla \phi(x)]+N,{~\rm where }\label{p1-DV-V}\\
%N&= (3-\alpha)t\|\dot x\|^2+(3-\alpha)t\|\dot \lambda\|^2- t\|Ax-b\|^2.\label{p1-DV-N}
%\end{eqnarray}
%
%\begin{eqnarray}\label{p1-DV-V}
%\dot V= 2t[\phi(x)-\phi(x^*)-(x-x^*)^{\top}\nabla \phi(x)]+N,
%\end{eqnarray}
%\begin{eqnarray}\label{p1-DV-N}
%where N= (3-\alpha)t\|\dot x\|^2+(3-\alpha)t\|\dot \lambda\|^2- t\|Ax-b\|^2.
%\end{eqnarray}

Because $\phi(\cdot)$ is convex, it is clear that $\phi(x)-\phi(x^*)-(x-x^*)^{\top}\nabla \phi(x)\leq 0$. It follows from \eqref{p1-DV-V} that
\begin{eqnarray}\label{DVleq0}
\dot V \leq (3-\alpha)t\|\dot x\|^2+(3-\alpha)t\|\dot \lambda\|^2- t\|Ax-b\|^2\leq 0,
\end{eqnarray}
which means that $V(t,x,\lambda,\dot x,\dot \lambda)$ is non-increasing on $[t_0,\infty)$. Hence, $$V_2(t,x,\dot x)\leq V(t,x,\lambda,\dot x,\dot \lambda)\leq V(t_0,x_0,\lambda_0,\dot x_0,\dot \lambda_0)\triangleq V_0.$$

Recall the definition of $V_2(\cdot)$. It follows that  $\|x+\beta t\dot x-x^*\|\leq \sqrt{\frac{V_0}{2}}$ and $\|x-x^*\|\leq \sqrt{\frac{V_0}{2(\alpha\beta-\beta-1)}}$ for all $t\geq t_0$. Clearly, $x(t)$ is bounded for all $t\geq t_0$. In addition, it follows from the  triangle inequality that $\|\beta t\dot x\|\leq \|x+\beta t\dot x-x^*\|+\|x^*-x\|\leq \sqrt{\frac{V_0}{2}}+\sqrt{\frac{V_0}{2(\alpha\beta-\beta-1)}}$ for all $t\geq t_0$. Hence, $\|t\dot x\|$ is bounded for all $t\geq t_0$. Similarly, one can prove that $\|\lambda\|$ and $\|t\dot \lambda\|$ are bounded for all $t\geq t_0$.
To sum up,  the trajectory of $(x(t),\lambda(t),t\dot x(t),t\dot \lambda(t))$ is bounded for $t\geq t_0$.  By \eqref{PD-Cent}, $(\ddot x(t),\ddot \lambda(t))$ is bounded for $t\geq t_0$.

(ii) Since $\dot V(t,x,\lambda,\dot x,\dot \lambda)\leq 0$, then $V(t,x(t),\lambda(t),\dot x(t),\dot \lambda(t))\leq V(t_0,x_0,\lambda_0,\dot x_0,\dot \lambda_0)\triangleq V_0$. It follows from the definition  of $V(\cdot)$  that $$t^2(L(x(t),\lambda^*)-L(x^*,\lambda^*))=V_1(t,x)\leq V(t,x(t),\lambda(t),\dot x(t),\dot \lambda(t))\leq V_0.$$ Hence,   $L(x(t),\lambda^*)-L(x^*,\lambda^*)=O(\frac{1}{t^2})$. Due to  the convexity of $\phi(\cdot)$  and \eqref{KKT1},  $$\phi(x)-\phi(x^*)+{\lambda^*}^\top A(x-x^*)\geq (x-x^*)^\top(\nabla \phi(x^*)+A^{\top}\lambda^*)=0.$$ Then it follows from  \eqref{Lag1} and $L(x(t),\lambda^*)-L(x^*,\lambda^*)=O(\frac{1}{t^2})$ that $\frac{1}{2}\|Ax(t)-b\|^2\leq L(x(t),\lambda^*)-L(x^*,\lambda^*)$ and $\|Ax(t)-b\|^2=O(\frac{1}{t^2})$. In addition, since we have proved the boundedness of $t\dot x(t)$ and $t\dot \lambda(t)$, it is clear that  $\|\dot x(t)\|=O(\frac{1}{t})$ and $\|\dot \lambda(t)\|=O(\frac{1}{t})$.

(iii) By part i), $(x(t),\lambda(t),t\dot x(t),t\dot \lambda(t))$ is bounded in a compact set for $t\geq t_0$. Because $\|Ax(t)-b\|^2\rightarrow 0$ as $t\rightarrow\infty$, $x(t)$ converges to the set of feasible points of problem \eqref{CP}.
Now, suppose, {\em ad absurdum}, that $x(t)$ does not converge to $S$ as $t\rightarrow\infty$. Recall that $x(t)$ stays in a compact set for all $t\geq t_0$. There exists an unbounded positive sequence $\{t_k\}\subset[t_0,\,\infty)$ such that $\lim_{k\rightarrow\infty}x(t_k)$ exists and $\mathrm{dist}(\lim_{k\rightarrow\infty}x(t_k),S)=\lim_{k\rightarrow\infty}\mathrm{dist}(x(t_k),S)\geq \epsilon>0$. Since $\|Ax(t)-b\|^2\rightarrow 0$ and  $L(x(t),\lambda^*)-L(x^*,\lambda^*)\rightarrow 0$ as $t\rightarrow\infty$, \begin{eqnarray*}
\lim_{k\rightarrow\infty}\phi(x(t_k))-\phi(x^*)&=&\phi(\lim_{k\rightarrow\infty}x(t_k))-\phi(x^*)\\
&=&L(\lim_{k\rightarrow\infty}x(t_k),\lambda^*)-L(x^*,\lambda^*)=0,
\end{eqnarray*} which contradicts $\lim_{k\rightarrow\infty}\mathrm{dist}(x(t_k),S)\geq \epsilon>0.$ Thus,  $x(t)$  converges to the solution set $S$ as $t\rightarrow\infty$.
%Furthermore, note that $0\leq V_1(t,x)\leq V(t,x,\lambda,\dot x,\dot \lambda)$. It follows that
%\begin{eqnarray}
%-\lambda^{*\top}A(x-x^*)-\frac{1}{2}\|Ax-b\|^2\leq \phi (x)-\phi(x^*)\leq \frac{1}{(t-t_0)^2}m_0-\lambda^{*\top}A(x-x^*)-\frac{1}{2}\|Ax-b\|^2.
%\end{eqnarray}
%Hence,
%\begin{eqnarray}
% |\phi (x)-\phi(x^*)|\leq \max\Big \{\frac{1}{(t-t_0)^2}m_0-\lambda^{*\top}A(x-x^*)-\frac{1}{2}\|Ax-b\|^2,\lambda^{*\top}A(x-x^*)+\frac{1}{2}\|Ax-b\|^2 \Big\}.
%\end{eqnarray}
%Since $\|Ax(t)-b\|^2=O(\frac{1}{t^2})$, $ |\phi (x)-\phi(x^*)|=O(\frac{1}{t^2}+\frac{1}{t})$.
%
%
%(iv)  Clearly,
%$$V(t,x(t),\lambda(t),\dot x(t),\dot \lambda(t))-m_0=\int_{t_0}^t \dot V(s,x(s),\lambda(s),\dot x(s),\dot \lambda(s)) \mathrm ds.$$
%Because $V(t,x(t),\lambda(t),\dot x(t),\dot \lambda(t))\geq 0$, $-\int_{t_0}^t \dot V(s,x(s),\lambda(s),\dot x(s),\dot \lambda(s)) \mathrm ds\leq m_0$. Inequalities \eqref{p1-ineq-int1}-\eqref{p1-ineq-int3} are ensured by  \eqref{p1-DV-V}.
\end{IEEEproof}

\begin{remark}
Main challenges  of proving Theorem \ref{thm-cent} are twofold. Firstly, one needs to find a suitable Lyapunov function to prove the convergence properties of  algorithm \eqref{PD-Cent} when $\alpha>3$. In the proof of Theorem  \ref{thm-cent}, the augmented  Lagrangian function and quadratic functions are elegantly combined to overcome this challenge.
Secondly, one needs to find suitable choices for $\beta$ when $\alpha >3$. In the proof of Theorem \ref{thm-cent}, $\beta=\frac{1}{2}$ is used to prove  \eqref{p1-DV-V}. In fact, by comparing  parameters in  \eqref{cp1-DVc1} and  \eqref{p1-temp}, one might  find that  $\beta=\frac{1}{2}$ is the only choice to obtain the convergence proof.
\hfill $\Diamond$
\end{remark}

One should note that the $O(\frac{1}{t^2})$ convergence rate is not tight  for strongly convex optimization.  Suppose  $\phi$ is a strongly convex or quadratic convex function. The primal-dual  dynamics has a linear rate, which is faster than $O(\frac{1}{t^2})$.

\subsection{ Analysis of $0< \alpha\leq 3$ and $\beta = \frac{3}{2\alpha}$}\label{Aleq3}

We proceed to discuss the case  with algorithm parameters $0<\alpha\leq 3$ and  $\beta= \frac{3}{2\alpha}$, and show in the following theorem  that the best rate of convergence is $O(t^{-\frac{2\alpha}{3}})$.
%\begin{subequations}\label{PD-Cents}
%\begin{align}
%\ddot x(t)\in & -\frac{\alpha}{t}\dot x(t)-\partial \phi(x(t))-A^{\top}(\lambda(t)+\beta t\dot \lambda(t))-A^{\top}(Ax(t)-b), \quad x(t_0)=x_0,\quad \dot x(t_0)=\dot x_0,\\
%\ddot \lambda(t)=& -\frac{\alpha}{t}\dot \lambda(t) + A(x(t)+\beta t\dot x(t))-b,\quad \lambda(0)=\lambda_0,\quad \dot \lambda(t_0)=\dot \lambda_0,
%\end{align}
%\end{subequations}

\begin{theorem}\label{thm-cents}
Denote $S$ as the set of solutions to  problem \eqref{CP} and assume $S\not=\emptyset$.   Let $(x(t),\lambda(t))$ be a trajectory of  algorithm \eqref{PD-Cent}, where $0<\alpha\leq 3$ and $\beta = \frac{3}{2\alpha}$.
\begin{itemize}
  \item [(i)] The trajectory of $( t^{\frac{\alpha}{3}-1}x(t), t^{\frac{\alpha}{3}-1}\lambda(t),t^{\frac{\alpha}{3}}\dot x(t),t^{\frac{\alpha}{3}}\dot \lambda(t))$ is bounded for $t\geq t_0$ and $0<\alpha<3$.
  \item [(ii)]The trajectory  $(x(t),\lambda(t),\dot x(t),\dot \lambda(t))$ satisfies the convergence properties
$L(x(t),\lambda^*)-L(x^*,\lambda^*)=O(t^{-\frac{2\alpha}{3}})$ and  $\|Ax(t)-b\|^2=O(t^{-\frac{2\alpha}{3}})$. If, in addition,   $0<\alpha< 3$, then $\|\dot x(t)\|=O(t^{-\frac{\alpha}{3}})$ and $\|\dot \lambda(t)\|=O(t^{-\frac{\alpha}{3}}).$
 \item  [(iii)] If, in addition, the trajectory of  $x(t)$ is bounded, then $x(t)$ converges to the solution set $S$ as $t\rightarrow\infty$.
\end{itemize}
\end{theorem}
\begin{IEEEproof}
To show that the rate of convergence is $O(t^{-\frac{2\alpha}{3}} )$, we define scalar $p\in(0,1]$ and functions $\theta:[t_0,\infty)\rightarrow\mathbb R_+$ and $\eta:[t_0,\infty)\rightarrow\mathbb R_+$.
Define the function
\begin{eqnarray}\label{Vcs}
V(t,x,\lambda,\dot x,\dot \lambda)=V_1(t,x)+V_2(t,x,\dot x)+V_3(t,\lambda,\dot \lambda)
\end{eqnarray}
such that
\begin{eqnarray*}
V_1&=& t^{2p}[L(x,\lambda^*)-L(x^*,\lambda^*)],\\
V_2 &=&0.5\|\theta(t)(x-x^*)+ t^p\dot x\|^2+\eta(t)/2\|x-x^*\|^2, \\
V_3 &=& 0.5\|\theta(t)(\lambda-\lambda^*)+t^p\dot \lambda\|^2+\eta(t)/2\|\lambda-\lambda^*\|^2,
\end{eqnarray*}
where $L(\cdot,\cdot)$ is defined in \eqref{Lag1} and $(x^*,\lambda^*)$ satisfies \eqref{KKT}.

(i) It follows from \eqref{KKT} and \eqref{PD-Cent} that derivatives of $V_i$'s along the trajectory of  algorithm \eqref{PD-Cent} satisfy
\begin{eqnarray}\label{DVc1s}
\dot V_1
&=& 2p t^{2p-1}[\phi(x)-\phi(x^*)+\lambda^{*\top}A(x-x^*)+\frac{1}{2}\|Ax-b\|^2]\notag\\
&&+ t^{2p}[\nabla \phi(x)+A^{\top}\lambda^*+A^{\top}(Ax-b)]^{\top}\dot x,
\end{eqnarray}
\begin{align}\label{DVc2s}
\dot V_2  =& \big (\theta(t)(x-x^*)+ t^p\dot x \big )^{\top}\Big (  \dot \theta(t)(x-x^*)+ (\theta(t)+ pt^{p-1})\dot x+t^p\ddot x \Big )\notag\\
&+\dot \eta(t)/2\|x-x^*\|^2+\eta(t)(x-x^*)^{\top}\dot x\notag\\
=& \big (\theta(t)(x-x^*)+ t^p\dot x \big )^{\top}\Big (  \dot \theta(t)(x-x^*)+ (\theta(t)+ (p-\alpha)t^{p-1})\dot x\notag\\
&+t^p(-\nabla \phi(x)-A^{\top}(\lambda+\beta t\dot \lambda)-A^{\top}(Ax-b)) \Big )\notag\\
&+\dot \eta(t)/2\|x-x^*\|^2+\eta(t)(x-x^*)^{\top}\dot x\notag\\
=& \Big(\theta(t)\dot \theta(t)+\frac{\dot \eta(t)}{2}\Big)(x-x^*)^2 + \Big (\theta^2(t)+ (p-\alpha)\theta(t)t^{p-1}+t^p\dot \theta(t)+\eta(t)\Big ) (x-x^*)^{\top}\dot x\notag\\
&-\theta(t)t^p(x-x^*)^{\top}(\nabla \phi(x)-\nabla \phi(x^*))-\theta(t)t^p(x-x^*)^{\top}A^{\top}(\lambda-\lambda^*)-\theta(t)t^p\|Ax-b\|^2\notag\\
&-\theta(t)\beta t^{p+1}(x-x^*)^{\top}A^{\top}\dot \lambda+t^p(\theta(t)+ (p-\alpha)t^{p-1})\dot x^2-t^{2p}(\nabla \phi(x)-\nabla \phi(x^*))^{\top}\dot x\notag\\
& -t^{2p}\dot x^{\top}A^{\top}(\lambda-\lambda^*)-\beta t^{2p+1}\dot x^{\top}A^{\top}\dot \lambda -t^{2p}\dot x^{\top}A^{\top}(Ax-b)),
\end{align}
\begin{align}\label{DVc3s}
\dot V_3  =& \big ( \theta(t)(\lambda-\lambda^*)+ t^p\dot \lambda \big )^{\top}\Big ( \dot \theta(t)(\lambda-\lambda^*)+ (\theta(t)+ pt^{p-1})\dot \lambda+t^p\ddot \lambda  \Big )\notag\\
&+\dot \eta(t)/2\|\lambda-\lambda^*\|^2+\eta(t)(\lambda-\lambda^*)^{\top}\dot \lambda\notag\\
=& \big ( \theta(t)(\lambda-\lambda^*)+ t^p\dot \lambda \big )^{\top}\Big ( \dot \theta(t)(\lambda-\lambda^*) + (\theta(t)+ (p-\alpha)t^{p-1})\dot \lambda+t^p( A(x+\beta t\dot x)-b) \Big )\notag\\
&+\dot \eta(t)/2\|\lambda-\lambda^*\|^2+\eta(t)(\lambda-\lambda^*)^{\top}\dot \lambda\notag\\
=&\Big(\theta(t) \dot \theta(t)+\frac{\dot \eta(t)}{2}\Big)(\lambda-\lambda^*)^2 + \Big (\theta^2(t)+ (p-\alpha)\theta(t)t^{p-1}+t^p\dot \theta(t)+\eta(t)\Big ) (\lambda-\lambda^*)^{\top}\dot \lambda\notag\\
&+\theta(t)t^p(\lambda-\lambda^*)^{\top}A(x-x^*)+\theta(t)\beta t^{p+1}(\lambda-\lambda^*)^{\top}A\dot x\notag\\
& +t^p(\theta(t)+ (p-\alpha)t^{p-1})\dot \lambda^2 +t^{2p}\dot \lambda^{\top}A(x-x^*)+\beta t^{2p+1}\dot \lambda^{\top}A\dot x.
\end{align}

Summing \eqref{DVc1s}-\eqref{DVc3s} and rearranging terms, it follows from \eqref{KKT} and \eqref{DVc1s}-\eqref{DVc3s} that
$\dot V(t,x,\lambda,\dot x,\dot \lambda) = \sum_{i=1}^5M_i,$
where
\begin{align*}
M_1=& 2p t^{2p-1}[\phi(x)-\phi(x^*)+\lambda^{*\top}A(x-x^*)+\frac{1}{2}\|Ax-b\|^2]\notag\\
&-\theta(t)t^p(x-x^*)^{\top}(\nabla \phi(x)-\nabla \phi(x^*))-\theta(t)t^p\|Ax-b\|^2,\notag\\
M_2=&\Big(\theta(t)\dot \theta(t)+\frac{\dot \eta(t)}{2}\Big)(x-x^*)^2+\Big(\theta(t) \dot \theta(t)+\frac{\dot \eta(t)}{2}\Big)(\lambda-\lambda^*)^2,\notag\\
M_3=& \Big (\theta^2(t)+ (p-\alpha)\theta(t)t^{p-1}+t^p\dot \theta(t)+\eta(t)\Big ) (x-x^*)^{\top}\dot x\notag\\
&+ \Big (\theta^2(t)+ (p-\alpha)\theta(t)t^{p-1}+t^p\dot \theta(t)+\eta(t)\Big ) (\lambda-\lambda^*)^{\top}\dot \lambda,\notag\\
M_4=&t^p\big (\theta(t)+ (p-\alpha)t^{p-1}\big )\dot x^2+t^p\big  (\theta(t)+(p-\alpha)t^{p-1}\big )\dot \lambda^2\notag,\\
M_5=&(t^{2p}-\theta(t)\beta t^{p+1})(x-x^*)^{\top}A^{\top}\dot\lambda +(\theta(t)\beta t^{p+1}-t^{2p})\dot x^{\top}A^{\top}(\lambda-\lambda^*).\notag
\end{align*}

Since $\nabla\phi(x^*) = -A^{\top}\lambda^*$ by \eqref{KKT1},   $M_1$ can be rewritten as:
\begin{align*}
&M_1= %2p t^{2p-1}[\phi(x)-\phi(x^*)+\lambda^{*\top}A(x-x^*)+\frac{1}{2}\|Ax-b\|^2]\notag\\
%&&-\theta(t)t^p(x-x^*)^{\top}(\partial \phi(x)-\partial \phi(x^*))-\theta(t)t^p\|Ax-b\|^2\\
%&=&
2p t^{2p-1}\Big [\phi(x)-\phi(x^*)-\frac{\theta(t)}{2p}t^{1-p}(x-x^*)^{\top}\nabla \phi(x)\\
&+\big (1-\frac{\theta(t)}{2p}t^{1-p}\big )\lambda^{*\top}A(x-x^*)+\big (\frac{1}{2}-\frac{\theta(t)}{2p}t^{1-p}\big )\|Ax-b\|^2\Big ].
\end{align*}
By considering $M_1$-$M_5$,  one can verify that $\dot V(t,x,\lambda,\dot x,\dot \lambda)\leq 0$ if
\begin{eqnarray}
\frac{\theta(t)}{2p}t^{1-p} &=& 1,\label{cond1}\\
\theta(t)\dot \theta(t)+\frac{\dot \eta(t)}{2}&\leq &0,\label{cond2}\\
\theta^2(t)+ (p-\alpha)\theta(t)t^{p-1}+t^p\dot \theta(t)+\eta(t)&=&0,\label{cond3}\\
\theta(t)+ (p-\alpha)t^{p-1}&\leq &0,\label{cond4}\\
t^{2p}-\theta(t)\beta t^{p+1}&=&0.\label{cond5}
\end{eqnarray}

Next, we seek feasible choices of  functions ($\theta(\cdot)$, $\eta(\cdot)$) and the parameter ($p$) that satisfy \eqref{cond1}-\eqref{cond5}.
Clearly, \eqref{cond1} implies that
\begin{eqnarray}\label{theta}
\theta(t) = 2p t^{p-1}.
\end{eqnarray}
Plugging \eqref{theta} in \eqref{cond3}-\eqref{cond5}, we have
\begin{eqnarray}
\eta(t)&=& (-2+\frac{\alpha+1}{2p})\theta^2(t),\label{eta1}\\
p &\leq & \frac{\alpha}{3}, {\rm~and ~} \beta  = \frac{1}{2p}. \label{p_alpha}
\end{eqnarray}
Plugging \eqref{theta} and \eqref{eta1} in \eqref{cond2} gives that
\begin{eqnarray}\label{ineq_alpha_p}
(-1+\frac{\alpha+1}{2p})(p-1)\leq 0.
\end{eqnarray}
Since $\eta(t)\geq 0$, \eqref{eta1} implies that $p\leq \frac{1+\alpha}{4}.$
 This combined with \eqref{p_alpha},   \eqref{ineq_alpha_p}, and   $p\in(0,1]$ proves that $0<p \leq  \min\{1, \frac{\alpha}{3}, \frac{1+\alpha}{4}\}=\frac{\alpha}{3}.$

Choose $p \triangleq \frac{\alpha}{3}$, $\theta(t) = \frac{2\alpha}{3} t^{\frac{\alpha}{3}-1}$, and $\eta(t)=\frac{2(3-\alpha)\alpha}{9}t^{\frac{2\alpha}{3}-2}$. Then $M_1$-$M_5$ can be simplified as
\begin{eqnarray*}
M_1 &=& \frac{2\alpha}{3} t^{\frac{2\alpha}{3}-1} [\phi(x)-\phi(x^*)-(x-x^*)^{\top}\nabla \phi(x) ]\leq 0,\\
M_2 &=&\frac{2}{27}\alpha(\alpha^2-9)t^{\frac{2\alpha}{3}-3}[(x-x^*)^2+(\lambda-\lambda^*)^2]\leq 0,
\end{eqnarray*}
$M_3 =M_4 =M_5 =0,$ and hence
\begin{eqnarray*}
\dot V & \leq & \frac{2\alpha}{3} t^{\frac{2\alpha}{3}-1} [\phi(x)-\phi(x^*)-(x-x^*)^{\top}\nabla \phi(x) ]\\
&&+\frac{2}{27}\alpha(\alpha^2-9)t^{\frac{2\alpha}{3}-3}[(x-x^*)^2+(\lambda-\lambda^*)^2]\leq 0.
\end{eqnarray*}

By plugging   $p = \frac{\alpha}{3}$, $\theta(t) = \frac{2\alpha}{3} t^{\frac{\alpha}{3}-1}$, and $\eta(t)=\frac{2(3-\alpha)\alpha}{9}t^{\frac{2\alpha}{3}-2}$ into $V_2(t,x,\dot x)$, we have
\begin{eqnarray*}
V_2 &=& {1\over 2}\Big{\|}\frac{2\alpha}{3} t^{\frac{\alpha}{3}-1}(x-x^*)+ t^{\frac{\alpha}{3}}\dot x\Big{\|}^2\notag+\frac{(3-\alpha)\alpha}{9}t^{\frac{2\alpha}{3}-2}\|x-x^*\|^2\\
%&=&\frac{\alpha^2+3\alpha}{9}t^{\frac{2\alpha}{3}-2}\|x-x^*\|^2+\frac{2\alpha}{3}t^{\frac{2\alpha}{3}-1}(x-x^*)^{\top}\dot x\notag\\
%&&+\frac{1}{2}t^{\frac{2\alpha}{3}}\|\dot x\|^2\\
&=&\Big{\|}\frac{\sqrt{\alpha^2+3\alpha}}{3}t^{\frac{\alpha}{3}-1}(x-x^*)+\sqrt{\frac{\alpha}{\alpha+3}}t^{\frac{\alpha}{3}}\dot x\Big{\|}^2+\frac{3-\alpha}{2(\alpha+3)}t^{\frac{2\alpha}{3}}\|\dot x\|^2.
\end{eqnarray*}

Since $V_1(t,x)\geq 0$  and function $V_2$ is positive definite with respect to  $(t^{\frac{\alpha}{3}-1}x,t^{\frac{\alpha}{3}-1}\lambda,t^{\frac{\alpha}{3}}\dot x,t^{\frac{\alpha}{3}}\dot \lambda)$ for all $t\geq t_0$ and $0<\alpha<3$. There exists a class $\mathcal K$ function $\kappa (\cdot)$ that $V(t,x,\lambda,\dot x,\dot \lambda)\geq \kappa (t^{\frac{\alpha}{3}-1}x,t^{\frac{\alpha}{3}-1}\lambda,t^{\frac{\alpha}{3}}\dot x,t^{\frac{\alpha}{3}}\dot \lambda)$.
It follows from $\dot V\leq 0$  that   the trajectory of $( t^{\frac{\alpha}{3}-1}x(t), t^{\frac{\alpha}{3}-1}\lambda(t),t^{\frac{\alpha}{3}}\dot x(t),t^{\frac{\alpha}{3}}\dot \lambda(t))$ is bounded for $t\geq t_0$ and $0<\alpha<3$.

(ii) Note that $V_2(t,x(t),\dot x(t)) \leq V(t,x(t),\lambda(t),\dot x(t),\dot \lambda(t))\leq V(t_0,x_0,\lambda_0,\dot x_0,\dot \lambda_0)$  for $t\geq t_0$. It follows that for $t\geq t_0$,
$\|\dot x(t)\|^2 \leq  \frac{2(\alpha+3)}{3-\alpha}\frac{1}{t^{\frac{2\alpha}{3}}}V_2(t,x(t),\dot x(t))
\leq \frac{2(\alpha+3)}{3-\alpha}\frac{1}{t^{\frac{2\alpha}{3}}}V(t_0,x_0,\lambda_0,\dot x_0,\dot \lambda_0).$
%\begin{eqnarray*}
%\|\dot x(t)\|^2& \leq & \frac{2(\alpha+3)}{3-\alpha}\frac{1}{t^{\frac{2\alpha}{3}}}V_2(t,x(t),\dot x(t))\\
%&\leq & \frac{2(\alpha+3)}{3-\alpha}\frac{1}{t^{\frac{2\alpha}{3}}}V(t_0,x_0,\lambda_0,\dot x_0,\dot \lambda_0).
%\end{eqnarray*}
Similarly, for $t\geq t_0$,
$\|\dot \lambda(t)\|^2 \leq \frac{2(\alpha+3)}{3-\alpha}\frac{1}{t^{\frac{2\alpha}{3}}}V(t_0,x_0,\lambda_0,\dot x_0,\dot \lambda_0).$
We have $\|\dot x(t)\|=O(t^{-\frac{\alpha}{3}})$ and $\|\dot \lambda(t)\|=O(t^{-\frac{\alpha}{3}}).$

As a result, $V(t,x(t),\lambda(t),\dot x(t),\dot \lambda(t))\leq V(t_0,x_0,\lambda_0,\dot x_0,\dot \lambda_0)$ for all $t\geq t_0.$
%$$V(t,x(t),\lambda(t),\dot x(t),\dot \lambda(t))-V(t_0,x_0,\lambda_0,\dot x_0,\dot \lambda_0)\leq \int_{s=t_0}^t \dot V(s,x(s),\lambda(s),\dot x(s),\dot \lambda(s)) \mathrm ds\leq 0.$$
Because $L(x(t),\lambda^*)-L(x^*,\lambda^*)=t^{-\frac{2\alpha}{3}} V_1(t,x(t))$ and
%\begin{eqnarray*}
$V_1(t,x(t)) \leq V(t,x(t),\lambda(t),\dot x(t),\dot \lambda(t))\leq V(t_0,x_0,\lambda_0,\dot x_0,\dot \lambda_0)$ for $t\geq t_0$,
%\end{eqnarray*}
we have  $L(x(t),\lambda^*)-L(x^*,\lambda^*)\leq t^{-\frac{2\alpha}{3}}V(t_0,x_0,\lambda_0,\dot x_0,\dot \lambda_0) $
%\begin{eqnarray*}
%L(x(t),\lambda^*)-L(x^*,\lambda^*)\leq t^{-\frac{2\alpha}{3}}V(t_0,x_0,\lambda_0,\dot x_0,\dot \lambda_0),
%\end{eqnarray*}
and   $L(x(t),\lambda^*)-L(x^*,\lambda^*)=O(t^{-\frac{2\alpha}{3}})$.
From  the definition of $L(\cdot,\cdot)$  in \eqref{Lag1}, it is clear that $\frac{1}{2}\|Ax(t)-b\|^2\leq L(x(t),\lambda^*)-L(x^*,\lambda^*)$ and $\|Ax(t)-b\|^2=O(t^{-\frac{2\alpha}{3}}).$

(iii) Suppose the trajectory of $x(t)$ is bounded. Because $\|Ax(t)-b\|^2\rightarrow 0$ as $t\rightarrow\infty$, $x(t)$ converges to the set of    feasible points of the problem \eqref{CP}. Then $L(x(t),\lambda^*)-L(x^*,\lambda^*)=O(t^{-\frac{2\alpha}{3}})$ implies that  $x(t)$  converges to the solution set $S$.
\end{IEEEproof}

\begin{remark}
%Theorem \ref{thm-cent} proves that both $L(x(t),\lambda^*)-L(x^*,\lambda^*)$ and $\|Ax(t)-b\|^2$ diminish at a convergence rate of $O(\frac{1}{t^2})$.
In Theorem \ref{thm-cents}, we show that the rate of convergence of $L(x(t),\lambda^*)-L(x^*,\lambda^*)= O(t^{-\frac{2\alpha}{3}})$ if $0< \alpha\leq 3$. However, the boundedness of $(x(t), \lambda(t), t\dot x(t), t\dot \lambda (t))$ is not guaranteed. In addition, combining the results of Theorem \ref{thm-cent}, $\alpha>3$ and $\beta=\frac{1}{2}$ are the optimal choice for parameters because they make algorithm \eqref{PD-Cent}  {converge with rate $O(\frac{1}{t^2})$.} %Hence, we can conclude from Theorems \ref{thm-cent} and \ref{thm-cents} that $\alpha>3$ is a better choice for obtaining a faster rate of convergence and the boundedness of algorithm trajectories.
\hfill $\Diamond$
\end{remark}

\begin{remark}
Different from the proof of Theorem \ref{thm-cent}, the construction of Lyapunov function $V$ in the proof of Theorem \ref{thm-cents} is more challenging because some time-varying gains $\theta(\cdot)$ and $\eta(\cdot)$ in the Lyapunov function need to be found. %Inspired by techniques in \cite{ACR:ESAIM:2017}, we
By analyzing necessary conditions  \eqref{cond1}-\eqref{cond5} for $\dot V\leq 0$, we find time-varying gains $\theta(\cdot)$ and $\eta(\cdot)$. If there is no dual variable  in algorithm  \eqref{PD-Cent}, the result is consistent with the result in \cite{ACR:ESAIM:2017} focusing the primal-based accelerated algorithm with $\alpha\leq 3$.
\hfill $\Diamond$
\end{remark}

The key of proving Theorems \ref{thm-cent} and \ref{thm-cents} is finding appropriate Lyapunov functions. The design of Lyapunov functions is partially inspired by the results for primal-based accelerated algorithms (see \cite{ACR:ESAIM:2017,Su2015A}). However, we have extended the algorithm design   and the analysis to primal-dual cases, which are a more general formulation. The duality gap function $L(x(t),\lambda^*)-L(x^*,\lambda^*)$ in the proof  is often used to analyze  primal-dual accelerated methods to show convergence orders (see \cite{Chambolle:2011:FPA:1968993.1969036,CLO:SIAM:2014,XZ:COA:2018}).
The obtained convergence rates, i.e., $O(\frac{1}{t^2} )$ for the case $\alpha>3$ and $O(t^{-\frac{2\alpha}{3}} )$ for the case $0<\alpha<3$,  are consistent to that of primal-based accelerated algorithms for unconstrained convex optimization problems \cite{ACR:ESAIM:2017,Su2015A}.

\subsection{Discussion on accelerated method}

In sections \ref{Ageq3} and \ref{Aleq3}, we have analyzed the primal-dual accelerated algorithm  for different choices of $\alpha$ and $\beta$. % and show that $\alpha>0$ and $\beta =\frac{1}{2}$ are better choices for algorithm \eqref{PD-Cent}.
The properties  are summarized in Table \ref{Tab1}, which shows that  $3$ is the critical value for $\alpha$.
Based on the analysis of lower bounds for convergence rates derived from constructed Lyapunov functions, we show that
convergence  rates of $L(x(t),\lambda^*)-L(x^*,\lambda^*)$, $\|Ax(t)-b\|^2$, $\|\dot x(t)\|$, and $\|\dot \lambda(t)\|$ {for the case $\alpha > 3$  are faster than those for the case $0<\alpha<3$}. In addition, if $\alpha > 3$, the right hand side of \eqref{PD-Cent} is bounded for any initial condition. Hence,  {$\alpha>3$} and $\beta =\frac{1}{2}$ are the best choices for parameters of \eqref{PD-Cent}.

\begin{table*}[htbp]\centering\caption{Algorithm performance for different $\alpha$}\label{Tab1}
	\label{table1}
	\begin{tabular}{|c|c|c|}
		\hline
		&$\alpha>3\quad \beta = \frac{1}{2}$ &$0<\alpha\leq 3\quad \beta = \frac{3}{2\alpha}$ \\
		\hline
		convergence rates for  $L(x(t),\lambda^*)-L(x^*,\lambda^*)$  and  $\|Ax(t)-b\|^2$   &$O(1/t^2)$& $O(1/t^{\frac{2\alpha}{3}} )$\\
		\hline
      convergence rates for $\|\dot x(t)\|$  and  $\|\dot \lambda(t)\|$ & $O(1/t)$ & $O( 1/t^{\frac{\alpha}{3}} )$\\
		\hline
% rates  of convergence for $ |\phi (x(t))-\phi(x^*)|$ & $O(\frac{1}{t^2}+\frac{1}{t})$ & $O(\frac{1}{t^{\frac{2\alpha}{3}}}+\frac{1}{t^\frac{\alpha}{3}})$ \\
%		\hline
   $ \|x(t)-x^*\|$ and $ \|\lambda(t)-\lambda^*\|$ & bounded  & not bounded \\
		\hline
%   {\color{blue}$ t^{(\frac{\alpha}{3}-1)}\|x(t)-x^*\|$  and $ t^{(\frac{\alpha}{3}-1)}\|\lambda(t)-\lambda^*\|$} & not bounded & bounded if $0<\alpha<3$ \\
%		\hline
  $ \|t\dot x(t)\|$ and $ \|t\dot \lambda(t)\|$ & bounded  & not bounded \\
		\hline
%   $ t^{\frac{\alpha}{3}}\|\dot x(t)\|$ and $t^{\frac{\alpha}{3}} \|\dot \lambda(t)\|$ & not bounded & bounded if $0<\alpha<3$ \\
%		\hline
   $ \|\ddot x(t)\|$ and $\|\ddot \lambda(t)\|$ & bounded & not bounded \\
		\hline
	\end{tabular}
\end{table*}

The proposed method is a continuous-time ordinary differential equation, which may be discretized to obtain discrete-time algorithms. %In fact, this is one motivation for research on continuous-time optimization algorithms.
However,  the discretization of the proposed method using explicit Euler scheme and modified symplectic scheme (see \cite{SDSJ:NIPS:2019}) has been  observed to be  numerically  instability by our simulation tests, which are omitted due to space limitations. The symplectic scheme \cite{SDSJ:NIPS:2019} may be a potential scheme to find a properly designed discrete-time counterpart for the proposed method.

%\begin{enumerate}
%  \item The rate of convergence for algorithm \eqref{AA}  is $O(\frac{1}{t^2})$, which is faster than that of algorithm \eqref{PD-Cents}.
%  \item Since the trajectory of $(x(t),y(t),t\dot x(t),t\dot y(t))$ generated by algorithm  \eqref{AA} is bounded, algorithm \eqref{AA} is well-defined with bounded $\ddot x(t)$ and $\ddot y(t)$. However,  $\ddot x(t)$ and $\ddot y(t)$ of algorithm \eqref{PD-Cents} may be unbounded as $t\rightarrow\infty$.
%\end{enumerate}

\section{Application to Network Optimization}
\label{DO}
In this section, we apply the proposed method \eqref{PD-Cent}  to  network optimization problems \eqref{S1} and \eqref{S2}, and design distributed primal-dual accelerated algorithms.

\subsection{Distributed Accelerated Algorithm for Scenario 1}
Consider  the  distributed primal-dual  accelerated algorithm
\begin{subequations}\label{DDPD1-Nest}
\begin{align}
\ddot x_i= & -\frac{\alpha_i}{t}\dot x_i-\nabla f_i(x_i)-\sum_{j=1}^n a_{i,j}(x_i-x_j)\notag\\
&-\sum_{j=1}^n a_{i,j}\big(\lambda_i+\frac{1}{2}t\dot \lambda_i-\lambda_j-\frac{1}{2}t\dot \lambda_j\big),\\
\ddot \lambda_i=& -\frac{\alpha_i}{t}\dot \lambda_i + \sum_{j=1}^n a_{i,j}\big (x_i+\frac{1}{2}t\dot x_i-x_j-\frac{1}{2}t\dot x_j\big ),
\end{align}
\end{subequations}
where  $t\geq t_0>0$, $x(t_0)=x_0,$ $\dot x(t_0)=\dot x_0,$ $\lambda(t_0)=\lambda_0,$ $\dot \lambda(t_0)=\dot \lambda_0,$ $\alpha_i>3$ is a parameter determined by agent $i\in\{1,\ldots,n\}$, and $a_{i,j}$ is the $(i,j)$th element of the adjacency matrix of graph $\mathcal G$.

Define $D_1 = \mathrm{diag}\{[\alpha_1,\ldots,\alpha_n]\}\otimes I_q$ and $L_{n\circ q} = L_n\otimes I_q$, where $L_n$ is the Laplacian matrix of $\mathcal G$.  Then the algorithm \eqref{DDPD1-Nest} has a compact formula as follows
\begin{subequations}\label{DPD1-Nest}
\begin{align}
\ddot x= & -\frac{1}{t}D_1\dot x-\nabla f(x)-L_{n\circ q}\big(\lambda+\frac{1}{2}t\dot \lambda\big)-L_{n\circ q}x, \\
\ddot \lambda=& -\frac{1}{t}D_1\dot \lambda + L_{n\circ q}\big ( x+\frac{1}{2}t\dot x\big).
\end{align}
\end{subequations}
By \cite[Lemma 3.1]{ZYH:2017}, we have the following result.
\begin{lemma}\label{lemma_4-1}
Let Assumption \ref{A1}  hold. Then $x^*\in\mathbb R^{nq}$ is a solution to   \eqref{S1} if and only if there exists $\lambda^*\in\mathbb R^{nq}$ such that
\begin{align}\label{S1KKT}
0_{nq}  = & \nabla f(x^*)+L_{n\circ q}\lambda^*  {\rm~and~}
0_{nq} =  L_{n\circ q}x^*.
\end{align}
\end{lemma}

Next, we present the main results of the algorithm \eqref{DPD1-Nest}.

\begin{theorem}\label{thm1}
Suppose  Assumption \ref{A1} holds. Let $(x(t),\lambda(t))$ be a trajectory of  algorithm \eqref{DPD1-Nest}.
Then
\begin{itemize}
  \item [(i)]  $(x(t),\lambda(t),t\dot x(t),t\dot \lambda(t))$ is bounded for $t\geq 0$;
  \item [(ii)]  $x(t)$ converges to the set of solutions to problem \eqref{S1}, and  $(x(t),\lambda(t),t\dot x(t),t\dot \lambda(t))$ satisfies the convergence properties
  $L_1(x(t),\lambda^*)-L_1(x^*,\lambda^*)=O(\frac{1}{t^2})$, $x^{\top}(t)L_{n\circ q}x(t)=O(\frac{1}{t^2})$, $\|\dot x(t)\|=O(\frac{1}{t})$, and $\|\dot \lambda(t)\|=O(\frac{1}{t})$.
\end{itemize}
\end{theorem}
\begin{IEEEproof}
Let $(x^*,\lambda^*)\in  \mathbb R^{nq}  \times \mathbb R^{nq} $ satisfy \eqref{S1KKT}. Define the augmented Lagrangian function of problem \eqref{S1} as
$L_1(x,\lambda) = f(x)+\lambda^{\top}L_{n\circ q}x.$
Define function
\begin{eqnarray}\label{V}
V(t,x,\lambda,\dot x,\dot \lambda)=V_1(t,x)+V_2(t,x,\dot x)+V_3(t,\lambda,\dot \lambda)
\end{eqnarray}
with $V_1 = \frac{1}{2}t^2\big [L_1(x,\lambda)-L_1(x^*,\lambda^*) +\frac{1}{2}x^{\top}L_{n\circ q}x\big ],$
\begin{align*}
V_2  =& \|x+\frac{t}{2}\dot x-x^*\|^2+\frac{1}{2}(x-x^*)^{\top}(D_1-3I_{nq})(x-x^*),\\
V_3  =& \|\lambda+\frac{t}{2}\dot \lambda-\lambda^*\|^2+\frac{1}{2}(\lambda-\lambda^*)^{\top}(D_1-3I_{nq})(\lambda-\lambda^*).
\end{align*}
%\begin{subequations}
%\begin{align}
%V_1 =&\frac{1}{2}t^2\big [L_1(x,\lambda)-L_1(x^*,\lambda^*) +\frac{1}{2}x^{\top}L_{n\circ q}x\big ]\notag\\
%=&\frac{1}{2}t^2\big [f(x)- f(x^*)+(x-x^*)^{\top}L_{n\circ q}\lambda^* +\frac{1}{2}x^{\top}L_{n\circ q}x \big ],\\
%V_2  =& \|x+\frac{t}{2}\dot x-x^*\|^2+\frac{1}{2}(x-x^*)^{\top}(D_1-3I_{nq})(x-x^*),\\
%V_3  =& \|\lambda+\frac{t}{2}\dot \lambda-\lambda^*\|^2+\frac{1}{2}(\lambda-\lambda^*)^{\top}(D_1-3I_{nq})(\lambda-\lambda^*).
%\end{align}
%\end{subequations}
Clearly, function $V$ is positive definite with respect to  $(x,\lambda,t\dot x,t\dot \lambda)$ for all $t\geq t_0>0$. One can prove the result using function $V$ as the proof of Theorem \ref{thm-cent}.
\end{IEEEproof}

\begin{remark}
The algorithm \eqref{DPD1-Nest} is an accelerated version of primal-dual algorithms in \cite{WE:2011,ZYH:2017}. Compared with \cite{WE:2011,ZYH:2017}, the rate of convergence for algorithm \eqref{DPD1-Nest} is accelerated to $O(\frac{1}{t^2})$. \hfill $\Diamond$
\end{remark}

\subsection{Distributed Optimization with Monotropic Constraint}

Consider the problem \eqref{S2}. We first decompose the coupled information in constraint \eqref{optimization_problem2}.
Define   $d \triangleq [d_1^{\top},\ldots,d_n^{\top}]^{\top}\in\mathbb R^{nm}$,  $z \triangleq [z_1^{\top},\ldots,z_n^{\top}]^{\top}\in\mathbb R^{nm}$, and $\overline{W} = \mathrm{diag}\{W_1,\ldots,W_n\}\in\mathbb{R}^{nm\times q}$.
Let $L_n$ be the Laplacian matrix of $\mathcal G$. Since $\mathcal G$ is connected and undirected by Assumption \ref{A1}, $\ker (L_{n})=\{v1_n: v\in\mathbb R\}$ and $\mathrm{range} (L_{n})=\{w\in\mathbb R^n:w^{\top}1_n=0\}$ form an orthogonal decomposition of $\mathbb R^{n}$ by the fundamental theorem of linear algebra \cite{Strang:1993}. Hence, $\sum_{i=1}^n W_i y_i =\sum_{i=1}^n d_i$ if and only if there exists $z\in \mathbb R^{nm}$ such that
$d-\overline{W} y-{L_{n\circ m}}{z}=0_{nm}$, where $L_{n\circ m} = L_{n}\otimes I_m$. It follows that problem \eqref{S2} is equivalent to the following problem
\begin{subequations}\label{S3}
\begin{align}
\min_{y\in\mathbb R^q,\ z\in\mathbb R^{nm}}&\,h(y), \quad h(y)= \sum_{i=1}^n h_i(y_i),\label{S3-1}\\
 \mathrm{s.t.}&\, d-\overline{W} y-{L_{n\circ m}}{z}=0_{nm}.\label{S3-2}
\end{align}
\end{subequations}

Let $ \lambda \triangleq [\lambda_1^{\top},\ldots,\lambda_n^{\top}]^{\top}\in\mathbb R^{nm}$.
We design the distributed primal-dual accelerated algorithm as
\begin{subequations}\label{PD-EMOD}
\begin{align}
\ddot y_i= & -\frac{1}{t}\alpha_i\dot y_i-\nabla h_i(y_i)+W_i^{\top}(\lambda_i+\frac{t}{2}\dot \lambda_i),\\
\ddot \lambda_i=& -\frac{1}{t}\alpha_i\dot \lambda_i +  d_i-W_i (y_i+\frac{t}{2}\dot y_i)-\sum_{j=1}^n a_{i,j} (\lambda_i-\lambda_j)-\sum_{j=1}^n a_{i,j} \big (z_i+\frac{t}{2}\dot z_i-z_j-\frac{t}{2}\dot z_j\big ),\\
\ddot z_i =& -\frac{1}{t}\alpha_i\dot z_i+\sum_{j=1}^n a_{i,j} \big (\lambda_i+\frac{t}{2}\dot \lambda_i-\lambda_j-\frac{t}{2}\dot \lambda_j\big ),
\end{align}
\end{subequations}
where $t\geq t_0>0$, $\alpha_i>3$, $y_i(t_0)=y_{i,0},$ $\dot y_i(t_0)=\dot y_{i,0}$, $\lambda_i(t_0)=\lambda_{i,0},$ $\dot \lambda_i(t_0)=\dot \lambda_{i,0},$ $z_i(t_0)=z_{i,0},$ $\dot z_i(t_0)=\dot z_{i,0},$ and $i\in\{1,\ldots,n\}$.

Define $D_2 \triangleq \mathrm{diag}\{\alpha_1 I_{q_1},\ldots,\alpha_n I_{q_n}\}$,  $D_3 \triangleq  \mathrm{diag}\{[\alpha_1,\ldots,\alpha_n]\}\otimes I_m$, and  the modified Lagrangian function $L_2:\mathbb R^{q}\times \mathbb R^{nm}\times \mathbb R^{nm}$ as
\begin{eqnarray}\label{Lagrng}
L_2(y,z,\lambda) = h(y)+\lambda^{\top}(d-\overline{W} y-{L}_{n\circ m}{z})-\frac{1}{2}\lambda^{\top}{L_{n\circ m}}{\lambda}.
\end{eqnarray}
The term $-\frac{1}{2}\lambda^{\top}{L_{n\circ m}}{\lambda}$ does not change the saddle points of the function due to \eqref{S2KKT}, and will help proving the convergence of $\lambda_i$'s. Hence, it is often used in the design of distributed algorithms for resource allocation problems (see \cite{ZYHX:SIAM:2018,YHL:Automatica:2016}).  Then the algorithm \eqref{PD-EMOD} is equivalent to the saddle point dynamics of \eqref{Lagrng} given by
\begin{subequations}\label{PD-EMO}
\begin{align}
\ddot y= & -\frac{1}{t}D_2\dot y-\nabla_{y} L_2\big (y,z,\lambda+\frac{t}{2}\dot \lambda\big), \\
\ddot \lambda=& -\frac{1}{t}D_3\dot \lambda +  \nabla_{\lambda} L_2\big (y+\frac{t}{2}\dot y,z+\frac{t}{2}\dot z,\lambda\big),\\
\ddot z =& -\frac{1}{t}D_3\dot z-\nabla_{z} L_2\big (y,z,\lambda+\frac{t}{2}\dot \lambda\big).
\end{align}
\end{subequations}
where $y(t_0)=y_0,$ $\dot y(t_0)=\dot y_0$, $\lambda(0)=\lambda_0,$ $\dot \lambda(0)=\dot \lambda_0,$ $z(t_0)=z_0$, and $\dot z(t_0)=\dot z_0$.

Similar to Lemma 4.3 of \cite{YHL:Automatica:2016,ZYHX:SIAM:2018}, the following result holds.
\begin{lemma}\label{lemma_4-2}
 Let  Assumption \ref{A1} hold. Then  $(y^*,z^*)\in\mathbb R^{nq}\times\mathbb R^{nm}$ is a solution to the problem  \eqref{S3} if and only if there exists $\lambda^*\in\mathbb R^{nq}$ such that
\begin{align}\label{S2KKT}
0_{q}  = & \nabla g(y^*)-\overline{W}^{\top}\lambda^*,~ 0_{nm} =   {L_{n\circ m}}{\lambda^*}, \notag \\
0_{nm} &=  d-\overline{W} y^*-{L_{n\circ m}}{z^*} .
\end{align}
\end{lemma}

The following theorem shows the convergence rate.
\begin{theorem}\label{conv-EMP}
Let  Assumption \ref{A1}  hold and $h(\cdot)$ be strictly convex. Suppose $(y(t),\lambda(t),z(t))$ is a trajectory of  algorithm \eqref{PD-EMO}. Then
\begin{itemize}
  \item [(i)] the trajectory of $(y(t),z(t),\lambda(t),t\dot y(t),t\dot z(t),t\dot \lambda(t))$ is bounded for $t\geq t_0>0$;
  \item [(ii)] the trajectory of $y(t)$ converges to the solution of problem \eqref{S2} and the trajectory satisfies
  $L_2(y(t),z^*,\lambda^*)-L_2(y^*,z^*,\lambda(t))=O(\frac{1}{t^2})$, $\lambda^{\top}(t)L_{n\circ m}\lambda(t)=O(\frac{1}{t^2})$, $\|\dot y(t)\|=O(\frac{1}{t})$,  $\|\dot \lambda(t)\|=O(\frac{1}{t})$, and $\|\dot z(t)\|=O(\frac{1}{t})$.
\end{itemize}
\end{theorem}
\begin{IEEEproof}
Let $(y^*,z^*,\lambda^*)$ satisfy \eqref{S2KKT}. Define $V(t,y,\lambda,z,\dot y,\dot \lambda,\dot z)=V_1(t,y)+V_2(t,y,\dot y)
 +V_3(t,\lambda,\dot \lambda)+V_4(t,z,\dot z)$
with $V_1= \frac{1}{2}t^2[ L_2(y,z^*,\lambda^*) -L_2(y^*,z^*,\lambda) ]$,
\begin{align*}
V_2 =& \|y+\frac{t}{2}\dot y-y^*\|^2+\frac{1}{2}(y-y^*)^{\top}(D_2-3I_{nm})(y-y^*),\\
V_3 =& \|\lambda+\frac{t}{2}\dot \lambda-\lambda^*\|^2+\frac{1}{2}(\lambda-\lambda^*)^{\top}(D_3-3I_{nm})(\lambda-\lambda^*),\\
V_4=& \|z+\frac{t}{2}\dot z-z^*\|^2+\frac{1}{2}(z-z^*)^{\top}(D_3-3I_{nm})(z-z^*),
\end{align*}
where $L_2(\cdot)$ is defined in \eqref{Lagrng}.  Then $V$ is positive definite with respect to  $(y,z,\lambda,t\dot y,t\dot z,t\dot \lambda)$ for all $t\geq t_0>0$.
Using $V$ as the Lyapunov function, one can prove this theorem similarly to that of Theorem \ref{thm-cent}. %Hence,  we omit  the proof.
\end{IEEEproof}

\begin{remark}
The algorithm \eqref{PD-EMO} is a modified version of the design in  \cite{ZYHX:SIAM:2018} by using the proposed accelerated method.
Compared with  results in  \cite{ZYHX:SIAM:2018}  which have $O(\frac{1}{t})$ convergence rate,  the algorithm \eqref{PD-EMO} has the $O(\frac{1}{t^2})$ convergence rate.\hfill $\Diamond$
\end{remark}
\subsection{Numerical Simulation}\label{simulation}
In this subsection, we present numerical simulations for the network optimization  problem \eqref{S2}.   Each local function is a log-sum-exp function
$h_i(y_i) = \rho \log \Big[ \sum_{j=1}^{m}\mathrm{exp}((c_{i,j}^{\top}y_i-b_{i,j})/\rho)   \Big],$
where $n=20$, $m=4$, $\rho=20$, $q = 2$, $d_0=[30, 50]^{\top}\in \mathbb R^{2}$, $ W_i\in\mathbb R^{2\times 2}$,  $c_{i,j}\in\mathbb R^{2}$ and $b_{i,j}\in\mathbb R$ are random vectors and scalars generated from a uniform distribution on the interval $[0,1]$.

\begin{figure}[htbp]
\centering
 \includegraphics[width=8 cm]{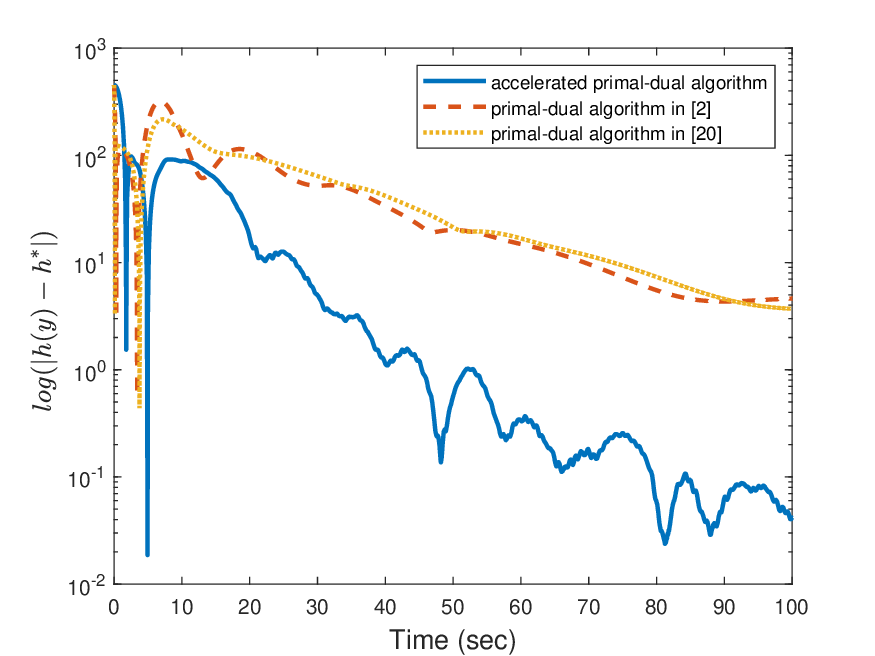}
  \caption{Trajectories of cost function $\mathrm{log}(|h(y)-h^*|)$
   %algorithm \eqref{PD-EMO} and   algorithms in \cite{ZYHX:SIAM:2018,MK:arxiv:2020}
  }
  \label{DCP1}
\end{figure}
\vspace{-0.2in}
\begin{figure}[htbp]
\centering
  \includegraphics[width=8 cm]{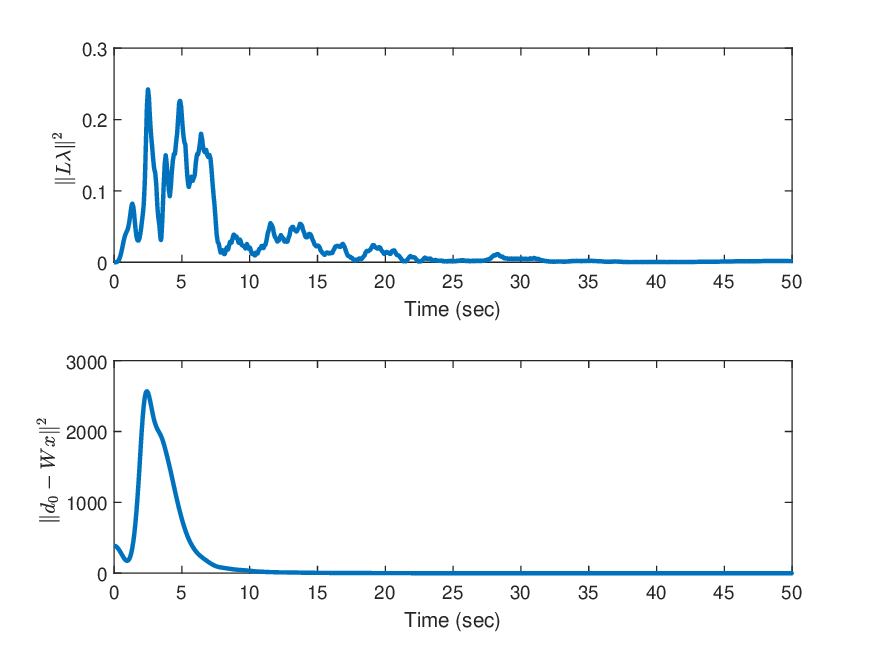}
  \caption{Trajectories of   $\|L_{n\circ m}\lambda\|^2$ and $\|d_0-Wy\|^2$ along algorithm \eqref{PD-EMO}}
  \label{DCP2}
\end{figure}

Simulation results of algorithm \eqref{PD-EMO}  are shown in Figs. \ref{DCP1} and \ref{DCP2}.  Because the $y$-axis of Fig. 1 is presented in log-scale, the slop of a trajectory shows the order of the convergence rate.  Fig. \ref{DCP1} shows that the accelerated design in this paper has a  faster performance than that of algorithms in  \cite{ZYHX:SIAM:2018} and \cite{MK:arxiv:2020}. In Fig. \ref{DCP2}, it is  shown that the dual variables will reach a  consensus and the constraint \eqref{optimization_problem2} is satisfied.

\section{Conclusion}\label{conclusion}
We have developed a primal-dual Nesterov's accelerated  method
for a class of convex optimization problems with affine equality constraints and applied the method to two types of network optimization problems. In particular, via a Lyapunov approach, we have analyzed critical choices of parameters in the algorithm design and proved that  the convergence rate of the Lagrangian function  is $O(\frac{1}{t^2})$, which is faster than the rate $O(\frac{1}{t})$ of standard primal-dual   methods. We further designed distributed accelerated primal-dual algorithms for optimization problems with consensus constraints and extended monotropic optimization problems using the proposed method.

Directions of future work include finding suitable discretization themes and designing convergent discrete-time counterparts  for the proposed continuous-time method. The analysis  of high-resolution ODE models and better measures of convergence rates to the accelerated primal-dual method is also of great importance. Moreover, the convergence to a minimizer rather than the set of solutions for primal-dual accelerated methods  is also worth investigating in future.

\bibliographystyle{IEEEtran}
\bibliography{Reference}

\end{document}